\documentclass[reqno, 12pt]{amsart}
 
\input cyracc.def
\font \tencyr = wncyr10
\def\cyr{\tencyr\cyracc}

\usepackage{amsmath, amsthm, amssymb}

\usepackage{color}

\newtheorem{thm}{Theorem}[section]
\newtheorem{cor}[thm]{Corollary}
\newtheorem{lem}[thm]{Lemma}

\theoremstyle{definition}
\newtheorem{rem}[thm]{Remark}

\theoremstyle{definition}

\theoremstyle{definition}
\newtheorem{assum}[thm]{Assumption}

\newcommand{\mysection}[1]{\section{#1}
\setcounter{equation}{0}}

\newcommand{\WO}[2]{\overset{\scriptscriptstyle0}{W}\,\!^{#1}_{#2}}

\makeatletter
\def\eint{\operatorname%
{\,\,\text{\bf--}\kern-.98em\DOTSI\intop\ilimits@\!\!}}
\makeatother

\def\coa{{\sf a}}
\def\cob{{\sf b}}
\def\coc{{\sf c}}

\def\bR{\mathbb{R}}

\title[Elliptic equations with measurable coefficients] {Elliptic
differential equations with measurable coefficients}
\author{Doyoon Kim \and N.V. Krylov}
\thanks{The work of the second author was partially supported by
NSF Grant DMS-0140405}

\address{127 Vincent Hall, University of Minnesota, Minneapolis, MN
55455} 
\email{dykim@math.umn.edu}

\address{127 Vincent Hall, University of Minnesota, Minneapolis, MN
55455} 
\email{krylov@math.umn.edu}


\begin{document}

\begin{abstract} We prove the unique solvability of second order elliptic
equations in non-divergence form in Sobolev spaces. The coefficients of
the second order terms are measurable in one variable and VMO in other
variables.  From this result, we obtain the weak uniqueness of the
martingale problem associated with the elliptic equations.
\end{abstract}

\maketitle

\mysection{Introduction}

We study the $L_p$-theory of the elliptic differential equation
\begin{equation}
                       		                             \label{11.13.1}
a^{jk}(x) u_{x^j x^k}(x) + b^{j}(x) u_{x^j}(x) + c(x)
u(x) = f(x)
\quad \text{in} \quad \bR^d,
\end{equation}
where $a^{jk}(x)$ are allowed to be only measurable with
respect to one coordinate, say $x^1 \in \bR$, where $x = (x^1, x') \in
\bR^d$, $x' \in
\bR^{d-1}$.

It is well known that if the coefficients $a^{jk}$ are only measurable, 
then there could not exist a {\em unique\/} solution to the above equation
even in a very generalized sense
(see \cite{Na, MR1684729}). We are interested in more regular solutions.
In 1967 Ural'tseva (see \cite{LU} or the original
paper
\cite{Ur})
constructed
an example of an equation in $\bR^{d}$ for $d\geq3$ with the
 coefficients
depending only on the first two coordinates for which
there is no unique solvability in  $W^{2}_{p}$ 
with $p\geq d$ (for any $d\geq3$ and $p\in (1,d)$   this was
 known before).

Thus in order to have the unique
solvability of the equation in $W^{2}_{p}$, we have to impose some
(regularity) conditions on the coefficients $a^{jk}$. The most classical
case is when
$a^{jk}$ are uniformly continuous. We can also have piecewise continuous
or VMO coefficients. For details, see \cite{
 MR1191890, MR1088476, Doyoon:article_half_cont:2004, Krylov_2005,
Lorenzi:article:1972, Softova}.

In this paper, we show that there exists a unique solution to the above
equation in $W_p^2 $, $p \in (2, \infty)$,  under the assumption
that $a^{jk}(x^1,x')$ are measurable in $x^1 \in \bR$ and VMO in $x' \in
\bR^{d-1}$. See Assumptions \ref{assum_01} and \ref{assum_02} below. 
 If the coefficients $a^{jk}$ are independent of $x' \in
\bR^{d-1}$ (more generally, uniformly continuous in $x' \in \bR^{d-1}$,
see Remark \ref{add_fact}), then the equation is uniquely solvable in
$W_2^2 $ as well.
In addition, we show that one can easily solve the  equation with
the Dirichlet,   Neumann, or oblique derivative boundary condition
in a half space, say $\bR^d_+ = \{(x^1, x') : x^1>0, x' \in \bR^{d-1}\}$,
using the results for equations in the whole space. 

The class of coefficients we are dealing with is considerably more
general than those previously known, as long as $p \in [2,\infty)$. It
actually contains almost all types of discontinuous coefficients that
have been investigated so far. For example, it contains the class of
piecewise continuous coefficients investigated in
\cite{Doyoon:article_half_cont:2004,Lorenzi:article_L_2:1972,
Lorenzi:article:1972}. It also contains VMO coefficients with which
elliptic equations were investigated in
\cite{MR1191890, MR1088476, Krylov_2005}.
  Also see the monograph \cite{Softova}, which treats elliptic and
parabolic equations with discontinuous coefficients including oblique
derivative problems with VMO coefficients. 
{\color{black} Although, we also slightly touch the oblique derivative
problem, we do not say anything about many important issues
of equations with VMO coefficients, which are discussed,
for instance, in \cite{So1}, \cite{So2}, \cite{PS}.

}

The highlight of our assumptions on the coefficients $a^{jk}$ would
be: no assumptions on the regularity of the coefficients
with respect to one variable as far as they are uniformly bounded and
elliptic. Having only measurable coefficients (as functions of $x^1 \in
\bR$), we obtain the
$L_2$-estimate for the equation by using the usual Fourier transforms.
Based upon this estimate, we establish the $L_p$-estimate, $p \in (2,
\infty)$, using the approach initiated by the second author of this paper
(for example, see \cite{Krylov_2005}). In this approach we make use of a
pointwise estimate of sharp functions of second order derivatives of the
solution.  As noted in \cite{Krylov_2005}, thanks to this method, we do
not need any integral representations of the solution nor commutators,
which were used, for example, in \cite{MR1191890, MR1088476}. 
Especially, we deal with VMO coefficients in a rather straightforward
manner. 

One good motivation to consider the above equation in the whole space is
to prove  weak uniqueness of stochastic processes associated with the
elliptic equation. As is shown in \cite{Krylov_2005, MR532498}, we can
say that   weak uniqueness of the processes holds true once we find a
unique solution of the elliptic equation in $W_p^2 $, 
$p \geq d$. More details are in \cite{Krylov_2005, MR532498}.

The paper is organized as follows. In Section \ref{main_results} we state
our main results.  The unique solvability of the equation in
$W_2^2 $  is investigated in  Section \ref{W_2^2_case}. In Section
\ref{pre_results}, we present some auxiliary results which are used in
Section \ref{W_p^2_case} where we finally prove the $W^{2}_p$-estimate,
$p
\in (2, \infty)$, for the equation.

{\color{black}
The authors are sincerely grateful to Hongjie Dong
who pointed out an omission in the first draft of
the article.
}
 
\mysection{Main results}\label{main_results}

We are considering the elliptic differential equation  \eqref{11.13.1}
where the coefficients $a^{jk}$, $b^{j}$, and $c$ 
 satisfy
the assumptions below.

\begin{assum}\label{assum_01} The coefficients $a^{jk}$, $b^{j}$, and $c$
are measurable functions defined on $\bR^d$, $a^{jk} = a^{kj}$. There
exist positive constants $\delta\in(0,1)$ and $K$ such that
$$ 
|b^{j}(x)|  \le K,
\qquad |c(x)| \leq K,
$$
$$
\delta |\vartheta|^2 \le \sum_{j,k=1}^{d} a^{jk}(x) \vartheta^j
\vartheta^k
\le\delta^{-1} |\vartheta|^2
$$ 
for any $x \in \bR^d$ and $\vartheta \in \bR^d$.
\end{assum}

To state another assumption on the coefficients, especially, $a =
(a^{jk})$, we introduce some notations. Let $B'_r(x') = \{ y' \in
\bR^{d-1}:  |x' - y'| < r \}$ and 
$Q_{r}(x)= Q_{r}(x^1,x') = (x^1-r, x^1+r) \times B'_r(x')$. Denote 
$$
\text{osc}_{x'}(a, Q_r(x)) = r^{-1} |B'_r|^{-2}
\int_{x^1-r}^{x^1+r}\int_{y', z' \in B'_r(x')} |a(t, y') - a(t, z') | \,
dy' \, dz' \, dt, 
$$
$$ 
a^{\#(x')}_R = \sup_{x \in \bR^d} \sup_{r \le R}
\text{osc}_{x'}(a, Q_r(x)),
$$ 
where $|B'_r|$ is the $d-1$-dimensional volume of
$B'_r(0)$. We  write  $a \in VMO_{x'}$ if 
$$
\lim_{R \to 0} a^{\#(x')}_R = 0.
$$ 
We see that $a \in VMO_{x'}$ if $a$ is independent of
$x'$.

\begin{assum}
                                           \label{assum_02} 

There is a continuous function $\omega(t)$
defined on $[0,\infty)$ such that $ \omega (0)=0$ and $a^{\#(x')}_R
\le
\omega(R)$ for all $R \in [0,\infty)$.
\end{assum} 

 \begin{rem}
It will be seen from our proofs that in Assumption
\ref{assum_02}  the requirement that $\omega(0)=0$
can be replaced with  $\omega(0) \leq (4 N_{1})^{-\nu (d+2)}$, where
$N_{1}=N_{1}(d,\delta,p)$ and $\nu = \nu(p)$ are the constants entering
\eqref{12.6.1}.

\end{rem}

As usual, we mean by $W_p^k(\bR^d)$, $k = 0, 1, \dots$, the Sobolev
spaces on $\bR^d$. Set $W_p^k = W_p^k(\bR^d)$, $L_p = L_p(\bR^d)$,
 and 
$$
Lu(x) = a^{jk}(x) u_{x^j x^k}(x) + b^{j}(x) u_{x^j}(x)
+ c(x) u(x).
$$ 
Here are our main results.

\begin{thm}\label{main_theorem} Let $p \in (2, \infty)$.
  Then there
exists a constant $\lambda_0$, depending only on $d$, $\delta$, $K$, $p$,
and the function $\omega$, such that, for any $\lambda 
>\lambda_0$ and
$f \in L_p$, there exists a unique $u \in W_p^2$ satisfying $Lu - \lambda
u = f$.

Furthermore, there is a constant $N$, depending only on $d$, $\delta$,
$K$, $p$, and the function $\omega$, such that, for any $\lambda \ge
\lambda_0$ and $u \in W_p^2$,
$$
\lambda \| u \|_{L_p} + \sqrt{\lambda} \| u_x \|_{L_p} + \| u_{xx}
\|_{L_p}
\le N \| Lu - \lambda u \|_{L_p}.
$$
\end{thm}

This theorem obviously covers the case in which the coefficients $a^{jk}$
are independent of $x' \in \bR^{d-1}$.  However, in that case we can
allow $p = 2$, which is detailed in the theorem below. Throughout the
paper, we write $N = N(d, \dots)$ if $N$ is a constant depending only on
$d,...$.
{\color{black}
The following theorem can be basically found in \cite{Ch}.
We give it a different proof that seems to be somewhat shorter
and more general.
}

\begin{thm}\label{W_2_case} Let the coefficients $a^{jk}$ be independent
of $x' \in \bR^{d-1}$.   Then there exists a constant $\lambda_0 =
\lambda_0(d, \delta, K)$ such that, for any $\lambda > \lambda_0$
  and
$f \in L_2$, there exists a unique $u \in W_2^2$ satisfying $Lu 
-\lambda u = f$.

In addition, there is a constant $N = N(d, \delta, K)$ such that, for any
$\lambda \ge \lambda_0$  and $u \in W_2^2$,
\begin{equation}                                    \label{L_2_esti}
\lambda \| u \|_{L_2} + \sqrt{\lambda} 
\| u_x \|_{L_2} + \| u_{xx}
\|_{L_2}
\le N \| Lu -\lambda u  \|_{L_2}.
\end{equation}
\end{thm}

\begin{rem}\label{add_fact} Theorem \ref{main_theorem} leads to the weak
uniqueness of solutions of stochastic differential equations associated
with the operator $L$. For details, see \cite{MR532498, Krylov_2005}. 
 Theorem \ref{W_2_case} clearly remains true under the assumption
that $a^{jk}(x^1,x')$ are uniformly continuous as functions of $x' \in
\bR^{d-1}$ uniformly in $x^1
\in \bR$. 
\end{rem}

 Three  more results
  deal with the equation $Lu - \lambda u = f$  in the half space 
$$
\bR^d_+ =\{x\in\bR^{d}:x^{1}>0\}.
$$
 Their proofs
 show the advantage of having the solvability in
$\bR^d$  of equations  whose coefficients are only  measurable in one
direction.  In what follows, we denote by $\WO{2}{p}(\bR^d_+)$  the
collection of all $u \in W_p^2(\bR^d_+)$ satisfying $u(0,x') \equiv 0$.

\begin{thm} 							\label{Dirichlet}
Let $p \in [2,\infty)$.  If $p = 2$, then suppose,
additionally, that the assumption in Theorem \ref{W_2_case} is satisfied.
Then there exists a constant $\lambda_0 = \lambda_0(d, \delta, K, p,
\omega){\color{black}
\geq0}$ such that, for any $\lambda > \lambda_0$ and
$f \in L_p(\bR^d_+)$, there exists a unique $u \in \WO{2}{p}(\bR^d_+)$
satisfying $Lu - \lambda u = f$.

Furthermore, there is a constant $N = N(d, \delta, K, p, \omega)$  such
that, for any $\lambda \ge
\lambda_0$ and $u \in \WO{2}{p}(\bR^d_+)$,
\begin{equation}                            \label{half_esti}
\lambda \| u \|_{L_p(\bR^d_+)} + \sqrt{\lambda} \| u_x \|_{L_p(\bR^d_+)} +
\| u_{xx}
\|_{L_p(\bR^d_+)}
\le N \| Lu - \lambda u \|_{L_p(\bR^d_+)}.
\end{equation}
\end{thm}

\begin{proof} We introduce a new operator 
$\hat{L} v = \hat{a}^{jk} v_{x^jx^k} + \hat{b}^j v_{x^j}
+ \hat{c} v$ the coefficients of which are as follows.
First we view the coefficients
$a^{jk}$, $b^{j}$, and $c$ as functions defined 
only on $\bR^d_+$. Then we
define $\hat{a}^{jk}$, $\hat{b}^{j}$, and $\hat{c}$ to be the odd or even
extensions of the 
original coefficients. Specifically, if $j=k=1$ or
$j,k
\in
\{ 2,
\dots, d \}$, then (even extension)
$$
\hat{a}^{jk}(x) =
\left\{\begin{array}{cl} a^{jk}(x^1,x') &\text{if} \,\,\, x^1 \ge 0 \\
a^{jk}(-x^1,x') &\text{if} \,\,\, x^1 < 0
\end{array}\right..
$$ If $j = 2, \dots, d$, then (odd extension)
$$
\hat{a}^{1j}(x) = \hat{a}^{j1}(x) =
\left\{\begin{array}{cl} a^{1j}(x^1,x') &\text{if} \,\,\, x^1 \ge 0 \\
-a^{1j}(-x^1,x') &\text{if} \,\,\, x^1 < 0
\end{array}\right..
$$ Similarly, the coefficient $\hat{b}^1(x)$ is the odd extension of
$b^1(x)$,  and the coefficients $\hat{b}^j(x)$, $j = 2, \dots, d$, and
$\hat{c}(x)$ are the even extensions of $b^j(x)$ and $c(x)$, respectively.
 
Now we notice that the coefficients of $\hat{L}$ satisfy Assumption
\ref{assum_01} and \ref{assum_02} with $2 \, \omega$. Then by Theorem
\ref{main_theorem} and \ref{W_2_case}, we can find a constant $\lambda_0 =
\lambda_0(d, \delta, K, p, \omega)$ such that, 
for any $\lambda > \lambda_0$, 
there exists a unique $u \in W_p^2$ satisfying $\hat{L} u -
\lambda u = \hat{f}$, where $\hat{f} \in L_p$ is the odd extension of $f
\in L_p(\bR^d_+)$. Obviously, $-u(-x^1,x') \in W_p^2$ also
satisfies the same equation, so by uniqueness we have $u(x^1,x') =
-u(-x^1,x')$. This implies that $u$, as a function defined on $\bR^d_+$,
is in the space
$\WO{2}{p}(\bR^d_+)$. 
Since $Lu - \lambda u = f$ in $\bR^d_+$, the function $u$ is a solution to the Dirichlet boundary problem.

To prove uniqueness and the estimate \eqref{half_esti},
we use the estimates in Theorem \ref{main_theorem} and \ref{W_2_case} and
the fact that the odd extension of an
element in $\WO{2}{p}(\bR^d_+)$ is in $W_p^2$.
The theorem is proved.
\end{proof}

In the same way, only this time taking the even extension of $f$, one gets
the solvability of the Neumann problem.
\begin{thm} Let $p \in [2,\infty)$.  If $p = 2$, then suppose,
additionally, that the assumption in Theorem \ref{W_2_case} is satisfied.
Then there exists a constant $\lambda_0 = 
\lambda_0(d, \delta, K, p, \omega){\color{black}
\geq0}$ such that, 
for any $\lambda > \lambda_0$ and
$f \in L_p(\bR^d_+)$, there exists a unique $u \in W^{2}_{p}(\bR^d_+)$
satisfying $Lu - \lambda u = f$ and $u_{x^{1}}=0$ on
$\partial\bR^{d}_{+}$.

Furthermore, there is a constant $N = N(d, \delta, K, p, \omega)$  such
that, for any $\lambda \ge
\lambda_0$ and $u \in W^{2}_{p}(\bR^d_+)$
satisfying  $u_{x^{1}}=0$ on
$\partial\bR^{d}_{+}$,
$$
\lambda \| u \|_{L_p(\bR^d_+)} + \sqrt{\lambda} \| u_x \|_{L_p(\bR^d_+)} +
\| u_{xx}
\|_{L_p(\bR^d_+)}
\le N \| Lu - \lambda u \|_{L_p(\bR^d_+)}.
$$
\end{thm}

While the Neumann problem is solved without any effort,
oblique derivative problems need some, still simple, manipulations.

Let $\ell$ be a constant vector field $\ell = 
(\ell^1, \cdots, \ell^d)$, where $\ell^1 > 0$.
Set $s = 1-1/p$ and recall that $g \in W_p^{s}(\bR^{d-1})$ if
$$
\|g\|_{W_p^{s}(\bR^{d-1})} 
= \|g\|_{L_p(\bR^{d-1})} + [g]_{s} < \infty,
$$
where
$$
[g]_{s}^p = \int_{\bR^{d-1}}\int_{\bR^{d-1}} 
\frac{|g(x') - g(y')|^p}{|x'-y'|^{d-1 + sp}} \, dx' \, dy'.
$$

\begin{thm}

                                             \label{theorem 12.3.1}

 Let $p \in [2,\infty)$.  If $p = 2$, then suppose,
additionally, that the assumption in Theorem \ref{W_2_case} is satisfied.
Then there exists a constant $\lambda_0 = 
\lambda_0(d, \delta, K, p, \omega, \ell){\color{black}
\geq0}$ such that, 
for any $\lambda >
\lambda_0$,
$f \in L_p(\bR^d_+)$, and $g \in W_p^{1-1/p}(\bR^{d-1})$, 
there exists a unique $u \in W^{2}_{p}(\bR^d_+)$
satisfying $Lu - \lambda u = f$ and $  \ell^j \, u_{x^{j}}=g$ on
$\partial\bR^{d}_{+}$.

Furthermore, there is a constant $N = N(d, \delta, K, p, \omega, \ell)$ such that, for any $\lambda \ge
\lambda_0$ and $u \in W^{2}_{p}(\bR^d_+)$,
\begin{multline}
                                                    \label{12.2.2}
\lambda \| u \|_{L_p(\bR^d_+)} + \sqrt{\lambda} \| u_x \|_{L_p(\bR^d_+)} +
\| u_{xx}
\|_{L_p(\bR^d_+)} \\
\le N \left(\| Lu - \lambda u \|_{L_p(\bR^d_+)} + (\lambda \vee
1)^{s/2} \|g\|_{L_p(\bR^{d-1})} + [g]_s \right),
\end{multline}
where  $\lambda \vee 1 = \max\{\lambda, 1\}$,  $s = 1-1/p$, 
and $g(x') =   \ell^j \, u_{x^{j}}(0,x')$.
\end{thm}

\begin{proof}

We can assume that $\ell^1 = 1$.
To introduce a new operator 
$$
\hat{L} v = \hat{a}^{jk} v_{x^jx^k} + \hat{b}^j v_{x^j} + \hat{c}
v,
$$
 we use a linear transformation 
$$
\varphi(x) = (-x^1, - 2 \ell' x^1
+ x')\quad (\ell' = (\ell^2, \dots, \ell^d)).
$$ 
Set
$$
\hat{a}^{jk}(x) =
\left\{\begin{array}{cl} a^{jk}(x) &\text{if} \,\,\, x^1 \ge 0 \\
\bar{a}^{jk}(x) &\text{if} \,\,\, x^1 < 0
\end{array}\right.,
$$ 
where 
$$
\bar{a}^{jk}(x) = \sum_{r,l=1}^d \varphi^j_{x^r} \varphi^k_{x^l} a^{rl}(\varphi(x)).
$$
Also set
$$
\hat{b}^j(x) = \left\{\begin{array}{cl} b^{j}(x) &\text{if} \,\,\, x^1 \ge 0 \\
\bar{b}^{j}(x) &\text{if} \,\,\, x^1 < 0
\end{array}\right., \quad
\hat{c}(x) = \left\{\begin{array}{cl} c(x) &\text{if} \,\,\, x^1 \ge 0 \\
\bar{c}(x) &\text{if} \,\,\, x^1 < 0
\end{array}\right.,
$$
where
$$
\bar{b}^j(x) = \sum_{r=1}^{d} \varphi^{j}_{x^r} b^{r}(\varphi(x)),
\qquad
\bar{c}(x) = c(\varphi(x)).
$$
Notice that the coefficients $\hat{a}^{jk}$ satisfy the uniform ellipticity condition with $N \delta$ in place of $\delta$, where $N$ depends only on $\ell$.
Also Assumption \ref{assum_02} is satisfied with $N \omega$ in place of
$\omega$, where $N$ depends only on $\ell$.

After this preparation we are ready to prove the first part
of the theorem.
Consider a differential equation 
\begin{equation}
                                                   \label{12.2.1}
\hat{L}w - \lambda w = \hat{f}_{\lambda}
\end{equation}
 in $\bR^d$, where
$\hat{f}_{\lambda}$ is defined as follows. 

One knows (see, for instance,  Theorem 2.9.1  of \cite{Tr}) that
for each $g\in W_p^{s}(\bR^{d-1})$ 
there is a function $v\in W^{2}_{p}(\bR^{d}_{+})$ such that
$v =0$ and $v_{x^{1}} =g(x')$ on $\partial\bR^{d}_{+}$ and, for a
constant
$N$ independent of $g$
\begin{equation}
													\label{extension_esti}
\|v\|_{W^{2}_{p}(\bR^{d}_{+})}\leq
N\|g\|_{W^{s}_{p}(\bR^{d-1})}.
\end{equation}
It follows by using dilations
that for any $g \in W_p^{s}(\bR^{d-1})$ and
$\lambda > 0$, we can find
$v
\in W_p^2(\bR^d_+)$ satisfying
$v  = 0$ and $v_{x^1}  = g $ on   $\partial\bR^{d}_{+}$, and
\begin{equation} 							       
                                                 \label{ext_esti}
\lambda \|v\|_{L_p(\bR^d_+)} + \sqrt{\lambda} \|v_x\|_{L_p(\bR^d_+)}
+ \|v_{xx}\|_{L_p(\bR^d_+)} \le N \left( \lambda^{s/2}\|g\|_{L_p(\bR^{d-1})} + [g]_{s} \right),
\end{equation}
where $N$ depends only on $d$ and $p$. We take this $v$
and  set
\begin{equation}							\label{extend_f}
\hat{f}_{\lambda}(x) = \left\{\begin{array}{cl}
f(x) - 2\hat{L}v + 2\lambda v &\text{if} \,\,\, x^1 > 0 \\
f(\varphi(x)) &\text{if} \,\,\, x^1 < 0
\end{array}\right..
\end{equation}
Using Theorem \ref{main_theorem} and \ref{W_2_case}, 
we find a unique solution $w \in W_p^2$ to \eqref{12.2.1}
for $\lambda > \lambda_0$, 
where $\lambda_0 = \lambda_0(d, \delta, K, p, \omega, \ell)$ is a constant corresponding to the operator $\hat{L}$.

Let $u^+$ be a function on $\bR^d_+$
defined by $u^+ = w + 2v$.
Also let $u^-$ be a function on 
$$
\bR^d_-= \{(x^1,x'): x^1 < 0, x' \in
\bR^{d-1}\}
$$ 
defined by $u^-(x) =
u^+(\varphi(x))$. We claim that $w = u^-$ in $\bR^d_-$. This
simple fact follows
from the uniqueness of solution to
the equation 
$\hat{L}w - \lambda w = \hat{f}_{\lambda}$ in $\bR^d_-$, proved in
Theorem \ref{Dirichlet}. Indeed, obviously,
  $w (0,x') = u^-(0,x')$ and it is also easy to check
  that $\varphi(\varphi(x))\equiv x$ and
$ \hat{L}u^- - \lambda u^- =
\hat{f}_{\lambda}$ in $\bR^d_-$.

Hence on $\partial\bR^{d}_{+}$
$$
w_{x^{1}} =u^{-}_{x^{1}} =(u^{+}(\varphi))_{x^{1}}
 =-u^{+}_{x^{1}} -2\sum_{j\geq2}\ell^{j}u^{+}_{x^{j}} .
$$
On the other hand, $w=u^{+}-2v$ on $\bR^{d}_{+}$ and
on $\partial\bR^{d}_{+}$
$$
w_{x^{1}} =u^{+}_{x^{1}} -2v_{x^{1}} 
=u^{+}_{x^{1}} -2g .
$$
It follows that on $\partial\bR^{d}_{+}$ it holds that
 $
\ell^{j}u^{+}_{x^{j}} =g 
 $
and since $ u^{+}  \in W^{2}_{p}(\bR^{d}_{+})$ we have
proved the existence of the desired solution.

To complete the proof, we now prove only \eqref{12.2.2}, which implies
uniqueness. Take a $u \in W_p^2(\bR^d_+)$ and set
  $g(x') =  \ell^j \, u_{x^j}(0,x')$.
Then for each 
$\lambda > \lambda_0 = \lambda_0(d, \delta, K, p, \omega,
\ell)$,  we find an extension $v \in W_p^2(\bR^d_+)$ of $g$
satisfying $v  = 0$, $v_{x^1} = g$ on $\partial\bR^{d}_{+}$,
and the estimate \eqref{ext_esti} or
\eqref{extension_esti}
depending on whether $\lambda \ge 1$ or $0 < \lambda < 1$.  
Define $w = u - 2v$ in $\bR^d_+$ and $w (x)  =
u(\varphi(x))$  in
$\bR^d_-$. 
Then $w(0+,x')=w(0-,x')$ 
 and  
$$
w_{x^{1}}(0+,x')=u_{x^{1}}(0,x')-2g(x')=-u_{x^{1}}(0,x')
-2\sum_{j\geq2} \ell^j u_{x^{j}}(0,x'),
$$
$$
w_{x^{1}}(0-,x')=-u_{x^{1}}(0,x')-2\sum_{j\geq2}
\ell^{j}u_{x^{j}}(0,x') =w_{x^{1}}(0+,x').
$$

It then follows that $w$  is a function in $W_p^2$ satisfying 
$\hat{L} w - \lambda w = \hat{f}_{\lambda}$, where $\hat{f}_{\lambda}$
is defined as in \eqref{extend_f} with $f:=Lu - \lambda u$. Hence by
Theorem \ref{main_theorem} and
\ref{W_2_case}, we have
$$
\lambda \| w \|_{L_p} + \sqrt{\lambda} \|w_x\|_{L_p} + \| w_{xx}
\|_{L_p}
\le N \|\hat{f}_{\lambda}\|_{L_p},
$$ where $N = N(d, \delta, K, p, \omega, \ell)$. 
 This, together with the estimates \eqref{extension_esti} and  
\eqref{ext_esti},  implies \eqref{12.2.2} 
for $\lambda>\lambda_{0}$. For $\lambda=\lambda_{0}$ we get
\eqref{12.2.2} by continuity.
\end{proof}

\begin{rem} Let $\ell(x') = (\ell^1(x'), \dots,
\ell^d(x'))$ be a bounded vector field defined on $\bR^{d-1}$ such
that $\ell(x') \in C^{1-1/p+\varepsilon}(\bR^{d-1})$, $\varepsilon >
0$, and $\ell^1(x') \ge \kappa > 0$. Then using the well-known
techniques -- freezing coefficients, partition of unity, and the
method of continuity, we can replace the constant vector field $\ell$
by $\ell(x')$ in the above theorem. Details can be found in
\cite{Softova}.
\end{rem}

\begin{rem}
A result similar to Theorem \ref{theorem 12.3.1}
holds if we replace the boundary condition
  $\ell^{j}u_{x^{j}}=g$ with $\ell^{j}u_{x^{j}}+\sigma u=g$,
where $\sigma$ is a constant. Indeed, again assuming
that $\ell^{1}=1$ it is easy to find
an infinitely differentiable bounded function $h(x^{1})$ 
having bounded derivatives and bounded away from zero
such that $h'(0)=-\sigma h(0)$. Then for $v=u/h$ we have
 $\ell^{j} v_{x^{j}}= g/h$  on $\partial\bR^{d}_{+}$ and
 $
Lu-\lambda u=h (\bar{L}v-\lambda v)$,
where
 $
\bar{L}\phi:=h^{-1}L(h \phi )
 $
is an elliptic operator satisfying our  hypotheses 
with  a slightly modified $K$.
\end{rem}

\mysection{Proof of Theorem 
\protect\ref{W_2_case}}\label{W_2^2_case}

Thanks
to the method of continuity and the denseness of
$C_0^{\infty}(\bR^d)$ in $W_2^2$, it suffices to prove the apriori
estimate \eqref{L_2_esti} for $u \in C_0^{\infty}(\bR^d)$
and   $a^{jk}$ that are sufficiently smooth. In
addition, on the account of possibly increasing $\lambda_{0}$ one sees
that it suffices  to prove   \eqref{L_2_esti}
for   $b\equiv0$,
$c\equiv0$, and $\lambda_{0}=0$.
In that case set
\begin{equation}
                                                  \label{10.31.1}
f=Lu-\lambda u.
\end{equation}

  For functions $\phi(x^{1},x')$ we denote by
$\tilde{\phi}(x^{1},\xi)$,
$\xi\in\bR^{d-1}$, its Fourier transform with respect to $x'$. 
  By taking the Fourier transforms 
 of both sides of \eqref{10.31.1}, we obtain
$$
\coa  \tilde{u}_{x^1x^1}  + \mathrm{i}   2   \cob 
\tilde{u}_{x^1}  - \coc  
\tilde{u}  = \tilde{f} ,
$$
\begin{equation}                                   \label{feq01}
 \tilde{u}_{x^1x^1}  + \mathrm{i}   2   \hat\cob 
\tilde{u}_{x^1}  - \hat\coc  
\tilde{u}  = \tilde{g} ,
\end{equation} where $\mathrm{i} = \sqrt{-1}$ and
$$
\coa(x^1) = a^{11}(x^1), 
\quad
\cob(x^1,\xi) = \sum_{j=2}^d a^{1j}(x^1) \xi^j,\quad\hat{\cob}
=\coa^{-1}\cob,
$$
$$
\coc(x^1,\xi) = \sum_{j,k=2}^d a^{jk}(x^1) \xi^j \xi^k+\lambda
,\quad\hat{\coc}
=\coa^{-1}\coc,\quad g=\coa^{-1}f.
$$

\begin{lem}\label{abcproperty} 
We have
$$
\delta \le \coa  = a^{11}  \le \delta^{-1},\quad 
|\cob(x^{1},\xi)|\leq \delta^{-1}|\xi|,
$$
\begin{equation} 
                                              \label{bdd a c}
\delta^{-1}(|\xi|^{2}+\lambda)\geq
 \coc(x^1,\xi) \geq \delta  |\xi|^2+\lambda ,
\end{equation} and
$$
\coa(x^1) \coc(x^1,\xi) - \cob^2(x^1,\xi) \ge \delta^2 (|\xi|^2
+\lambda).
$$
\end{lem}

\begin{proof} We prove only the last inequality. From Assumption
\ref{assum_01}, we have
$$
\delta ( t^2 + |\xi|^2)
\le \coa(x^1) \, t^2 + 2 \, \cob(x^1,\xi) \, t + \coc(x^1,\xi)-\lambda.
$$ 
for all $t \in \bR$ and $\xi \in \bR^{d-1}$.
 In particular, 
$$
\left( \coa(x^1) - \delta \right) t^2 + 2 \, \cob(x^1,\xi) \, t +
\coc(x^1,\xi) - \delta |\xi|^2-\lambda\ge 0.
$$ 
This implies that
$$
\cob^2(x^1,\xi)  - \left( \coa(x^1) - \delta \right) \left( \coc(x^1,\xi)
- \delta |\xi|^2 -\lambda\right) \le 0.
$$ 
 From this and \eqref{bdd a c} the result follows.
\end{proof}

\begin{lem}                          \label{midimp01}
 For any
$\xi\in\bR^{d}$
\begin{equation}
                                                  \label{10.28.1}
  (|\xi|^2 + \lambda)   \int_{\bR}    |\tilde{u}_{x^{1}}|^{2}\,dx^{1} + 
(|\xi|^4+\lambda|\xi|^2+\lambda^{2})
\int_{\bR}   |\tilde{u}|^{2}\,dx^{1}
\le N(\delta ) \int_{\bR} |\tilde{f}|^2 \, dx^1,
\end{equation}
\begin{equation}
                                                  \label{10.28.2}
 \int_{\bR}    |\tilde{u}_{x^{1}x^{1}}|^{2}\,dx^{1}
\le N(\delta ) \int_{\bR} |\tilde{f}|^2 \, dx^1.
\end{equation}
\end{lem}

\begin{proof} Estimate \eqref{10.28.2} is a direct consequence of
equation \eqref{feq01} (allowing one to express $\tilde{u}_{x^{1}x^{1}}$
through $\tilde{f}$, $\tilde{u}_{x^{1}}$, and
$\tilde{u}$), \eqref{bdd a c}, and \eqref{10.28.1}.

While proving \eqref{10.28.1} we
define a function $\phi(x^1, \xi)$ by
$\phi(0,\xi)=0$ and $ \phi_{x^{1}}=\hat\cob$  and set $\rho  = \tilde{u}
e^{\mathrm{i} \phi }$.  Then from \eqref{feq01} we see that
$$
\rho_{x^1x^1}  +  (
 \hat\cob ^2   - \mathrm{i} \, \phi_{x^1x^1}   
  -  \hat\coc   
 ) \rho    = \tilde{g}  e^{\mathrm{i} \phi }.
$$ 
Multiply both sides by $\bar{\rho}$ and integrate the result with
respect to $x^{1}$. Integrating by parts shows that
$$ 
- \int_{\bR} |\rho_{x^1}|^2 \, dx^1 + \int_{\bR} (
 \hat\cob ^2   - \mathrm{i} \, \phi_{x^1x^1}   
  - \hat\coc  
 ) |\tilde{u}|^2 \, dx^1  = \int_{\bR}
 \tilde{g}\bar{\tilde{u}}\, dx^1.
$$ 
Taking the real parts of both sides and multiplying by
$|\xi|^2+\lambda$, we have
$$
\int_{\bR}(|\xi|^2+\lambda) |\rho_{x^1}|^2 \, dx^1 +
\int_{\bR}  (\hat \coc - \hat\cob^2) (|\xi|^2+\lambda) |\tilde{u}|^2
\, dx^1  
$$
$$
= -\int_{\bR}  (|\xi|^2+\lambda) 
 \Re (\tilde{g}\bar{\tilde{u}}) \, dx^1.
$$
 Note that for any $\varepsilon>0$
$$
 - (|\xi|^2+\lambda) \Re (\tilde{g}\bar{\tilde{u}})
\le  \varepsilon  (|\xi|^2+\lambda)^2 |\tilde{u}|^2 + 
  \varepsilon^{-1} |\tilde{g}|^2.
$$ 
 From this and Lemma \ref{abcproperty} we obtain
$$
\int_{\bR}  (|\xi|^2 + \lambda) |\rho_{x^1}|^2 \, dx^1 +
\int_{\bR} ( \delta^4 -\varepsilon) 
 (|\xi|^2+\lambda)^{2} |\tilde{u}|^2 \, dx^1 
\le  \varepsilon^{-1} \int_{\bR} |\tilde{g}|^2 \, dx^1.
$$
 By choosing an appropriate $\varepsilon > 0$  (e.g.
$\varepsilon = \delta^4/2 $), we arrive at
$$
\int_{\bR}  (|\xi|^2 + \lambda)  |\rho_{x^1}|^2 \, dx^1 + \int_{\bR}  
(|\xi|^4+\lambda|\xi|^2+\lambda^{2})|\tilde{u}|^2 \, dx^1 
\le N(\delta ) \int_{\bR} |\tilde{f}|^2 \, dx^1.
$$

It only remains to observe that in light of \eqref{bdd a c}
$$
 |\tilde{u}_{x^{1}} | = |\rho_{x^{1}}-\mathrm{i}
\cob\coa^{-1}\tilde{u}e^{i\phi}| \leq|\rho_{x^{1}}|+
N(\delta)|\xi||\tilde{u}|.
$$~\end{proof}

 Now we can finish the proof  of Theorem ~\ref{W_2_case}.
As we pointed out in the beginning of the section
we only need to prove \eqref{L_2_esti} for
$u\in C^{\infty}_{0}(\bR^{d})$, 
smooth $a^{ij}$, $b\equiv0$, $c\equiv0$, and 
$\lambda_{0}=0$.

 In that
case 
it suffices to add  \eqref{10.28.1} and \eqref{10.28.2}, integrate over
$\bR^{d-1}$ and use  Parseval's identity. 
 The theorem is proved.

\begin{rem}\label{general L_2} We  have just proved that 
if $b^{j} = c = 0$, then  
$$
\lambda \| u \|_{L_2} + \sqrt{\lambda} \| u_x \|_{L_2} + \| u_{xx}
\|_{L_2}
\le N \| L u - \lambda u \|_{L_2}
$$ 
for $u \in W_2^2$ and $\lambda \ge 0$, where $N$ depends
only on  $\delta$. 
\end{rem}

\mysection{Auxiliary results}\label{pre_results}

Here we state and prove a series of observations which are needed in the
proof of Theorem \ref{main_theorem}. First we introduce some notation.
 As usual, we set $B_r(x_0) = \{ x \in \bR^d : |x-x_0| < r \}$ and $B_r =
B_r(0)$. By $|B_r|$ we mean the $d$-dimensional volume of
$B_r$. We denote by $|u|_0$ the supremum of $u$ over $\bR^d$.
 
 Throughout this section, we assume that 
$$ 
L u(x) = L_0 u(x) = a^{jk}(x^1) u_{x^j x^k}(x).
$$

Our first auxiliary result is the following.
\begin{lem}                                   \label{L_2 estimate}
 There exists  $N = N(d,\delta)$  
such that, for any $u \in W_2^2(B_{R})$ with
$u|_{\partial B_{R}} = 0$, we have
\begin{equation}                               \label{main L_2}
R^{2}\int_{B_R} |u_{x}|^2 \, dx +
\int_{B_R} |u|^2 \, dx 
\le N \, R^{4 }
 \int_{B_{R}} |Lu|^{2 } \, dx  .
\end{equation}
\end{lem}

\begin{proof}
Assume that \eqref{main L_2} is true when $R = 1$. For a
given $u \in W_2^2(B_R)$ with $u|_{\partial B_R} = 0$, we set 
$$ 
L_R = a^{jk}(Rx )  \frac{\partial^2}{\partial x^j
\partial x^k} 
\quad \text{and} \quad v(x) = R^{-2} u(Rx ).
$$ 
Then $v \in W_2^2(B_1)$ and $L_R v(x) = (Lu)(Rx )$ in
$B_1$. Since $L_R$ satisfies the same ellipticity condition as  
$L$ does, we have
$$
\int_{B_R} |u|^2 \, dx = R^{d+4} \int_{B_1} |v|^2 \, dx 
$$
$$
\le N R^{d+4}   \int_{B_1} |L_R v|^{2} \, dx
 = N R^{4  }  \int_{B_R}
|Lu|^{2} \, dx  .
$$
Also
$$
\int_{B_R} |u_{x}|^2 \, dx = R^{d+2} \int_{B_1} | v_{x} |^2 \, dx 
$$
$$
\le N R^{d+2}   \int_{B_1} |L_R v|^{2} \, dx
 = N R^{2}  \int_{B_R}
|Lu|^{2} \, dx  .
$$
This shows that we need only  prove the lemma for
$R=1$.

In that case we can divide $L$ by $a^{11}$
and may assume that $a^{11}\equiv1$.
Then we   integrate $uLu$ over $B_{1}$ using integration
by parts to find
$$
\delta\int_{B_{1}}|u_{x}|^{2}\,dx
\leq\int_{B_{1}}a^{jk}u_{x^{j}}u_{x^{k}}\,dx
=- \int_{B_{1}}uLu\,dx
$$
$$
\leq\big(\int_{B_{1}}u^{2}\,dx\big)^{1/2}
\big(\int_{B_{1}}(Lu)^{2}\,dx\big)^{1/2}.
$$
We estimate the integral of $u^{2}$ through that of $|u_{x}|^{2}$
by using Poincar\'e's inequality and obtain the needed
estimate for $u_{x}$. This is the only estimate we
need to prove since $u$ is estimated by $u_{x}$ again owing to
Poincar\'e's inequality.

\end{proof}

The following lemma is almost identical to a theorem in
\cite{Kr67}.  For  completeness, we present here a proof.

\begin{lem}\label{estimate K} Let $0 < r < R$. There exists 
 $N = N(d,\delta)$  
such that, for $w \in W_2^2(B_R)$, 
$$
\| w \|_{W_2^2(B_r)} \le N \left( \| Lw - w\|_{L_2(B_R)} + (R-r)^{-2}
\|w\|_{L_2(B_R)} \right).
$$
\end{lem}

\begin{proof} Let 
$$ 
R_0 = r,
\quad R_m = r + (R-r) \sum_{k=1}^m \frac{1}{2^k}, 
\quad m = 1, 2, \dots,
$$ 
$$ 
B_m = \{ x \in \bR^d : |x| \le R_m \},
\quad m = 0, 1, \dots.
$$ 
Also let
$\zeta_m \in C_0^{\infty}(\bR^d)$ such that $\zeta_m(x) = 1$ in $B_m$,
$\zeta_m(x) = 0$ outside of $B_{m+1}$, and
$$ 
| (\zeta_m)_x |_0 \le N \frac{2^{m+1}}{(R-r)},
\quad | (\zeta_m)_{xx} |_0 \le N \frac{2^{2m+2}}{(R-r)^2},
$$ 
where $N$ depends only on $d$. To construct them take an
infinitely differentiable function $g(t)$, $t \in (-\infty, \infty)$,
such that $g(t) = 1$ for $t \le 1$, $g(t) = 0$ for $t \ge 2$, and $0 \le
g \le 1$. After this define 
$$
\zeta_m(x) = g( 2^{m+1}(R-r)^{-1} (|x| - R_m) + 1 ).
$$

Now we make use of the $L_2$-estimate of $\zeta_m w$, which is from
Remark~\ref{general L_2}, as follows.
\begin{multline}\label{apply L2}
\|w\|_{W_2^2(B_m)}
\le \| \zeta_m w \|_{W_2^2}
\le N \| (L-1) \zeta_m w \|_{L_2} \\
\le N \| (L-1) w \|_{L_2(B_R)} + N \frac{2^{m+1}}{R-r} \| w_x
\|_{L_2(B_{m+1})} + N \frac{2^{2m+2}}{(R-r)^2} \| w \|_{L_2(B_R)},
\end{multline} where $N$ depends only on  $d$ and $\delta$. 
By interpolation inequalities
$$
\| w_x \|_{L_2(B_{m+1})}
\le \varepsilon \|w_{xx}\|_{L_2(B_{m+1})} + N \varepsilon^{-1}
\|w\|_{L_2(B_{m+1})},
$$ 
where $\varepsilon > 0$, and $N$ depends only on $d$ (by
using a dilation argument we can take a constant $N$ which does not
depend on the radius of $B_{m+1}$). Thus the right hand side of the
inequality \eqref{apply L2} is not greater than
\begin{equation}\label{app inter} N \| (L-1) w \|_{L_2(B_R)} +
\varepsilon \|w_{xx}\|_{L_2(B_{m+1})} + N \varepsilon^{-1}
\frac{2^{2m+2}}{(R-r)^2} \|w\|_{L_2(B_R)},
\end{equation} where $0 < \varepsilon < 1$ and $N$ depends only on 
 $d$ and $\delta$. Set 
$$
\mathcal{A}_m := \|w\|_{W_2^2(B_m)}, 
\quad
\mathcal{B} := \| (L-1) w \|_{L_2(B_R)},
\quad \text{and} 
\quad
\mathcal{C} := \|w\|_{L_2(B_R)}.
$$ 
Then from \eqref{apply L2} and \eqref{app inter}, we have
$$
\varepsilon^m \mathcal{A}_m \le N \varepsilon^m \mathcal{B} +
\varepsilon^{m+1} \mathcal{A}_{m+1} + N \varepsilon^{m-1}
\frac{2^{2m+2}}{(R-r)^2} \mathcal{C}.
$$ 
Choose an $\varepsilon$ such that $0 <  4 \varepsilon <
1$, and notice that
$\mathcal{A}_m \le \|w\|_{W_2^2(B_R)}$.  Then we   have
$$
\sum_{m=0}^{\infty}\varepsilon^m \mathcal{A}_m \le N \mathcal{B}
\sum_{m=0}^{\infty} \varepsilon^m + \sum_{m=0}^{\infty} \varepsilon^{m+1}
\mathcal{A}_{m+1} + N \frac{\varepsilon^2}{(R-r)^2} \mathcal{C}
\sum_{m=0}^{\infty} (4 \varepsilon)^{m+1}.
$$ 
This clearly finishes the proof.
\end{proof}

\begin{rem}\label{lambda 0} Using the dilation argument as in the proof
of Lemma~\ref{L_2 estimate}, we   have
\begin{multline*}
\lambda \| w \|_{L_2(B_r)} + \sqrt{\lambda}\| w_x \|_{L_2(B_r)} + \|
w_{xx} \|_{L_2(B_r)} \\
\le N \left( \| Lw - \lambda w\|_{L_2(B_R)} + (R-r)^{-2} \|w\|_{L_2(B_R)}
\right)
\end{multline*} for any $\lambda > 0$, where $N$ depends only on 
 $d$ and $\delta$. In particular, by letting $\lambda \to 0$, we have
\begin{equation}\label{case lambda0}
\| w_{xx} \|_{L_2(B_r)} \le N \left( \| Lw \|_{L_2(B_R)} + (R-r)^{-2}
\|w\|_{L_2(B_R)} \right).
\end{equation}
\end{rem}

In the next few lemmas, we investigate some properties of a solution $h$
of the equation $Lh = 0$. Recall that the coefficients $a^{jk}$ of the
operator $L$ do not depend on $x' \in \bR^{d-1}$.

\begin{lem}\label{conse esti} Let $\gamma = (\gamma^1, \dots, \gamma^d)$
be a multi-index such that
$\gamma^1 = 0, 1, 2$. Also let $0 < r < R \le 4$. If $h$ is a
sufficiently smooth function defined on $B_4$ such that $Lh = 0$ in $B_4$,
then we have 
$$
\int_{B_r} |D^{\gamma} h|^2 \, dx
\le N \int_{B_R} |h|^2 \, dx,
$$ 
where  $N = N(d, \delta, \gamma, R, r)$.
\end{lem}

\begin{proof} Set $\gamma' = (0, \gamma^2, \dots, \gamma^d)$ and notice
that
$$ 
L (D^{\gamma'}h) = 0, 
\quad \text{that is}, \quad (L - 1)D^{\gamma'}h = - D^{\gamma'}h 
\quad \text{in} \quad B_4.
$$ 
Then by Lemma \ref{estimate K}
$$
\| D^{\gamma}h \|_{L_2(B_r)}
\le N \left( \|D^{\gamma'}h\|_{L_2(B_{r_1})} + (r_1 - r)^{-2}
\|D^{\gamma'}h\|_{L_2(B_{r_1})} \right),
$$ 
where $r < r_1 < R$. If $|\gamma'| = 0$, then we are done.
Otherwise, we can consider a multi-index $\gamma''$
having at least one component less by one than the corresponding
component of $\gamma'$.
Then,   $L(D^{\gamma''}h) = 0$ and
$$
\|D^{\gamma'} h\|_{L_2(B_{r_1})}
\le N \left( \|D^{\gamma''}h\|_{L_2(B_{r_2})} + (r_2 - r_1)^{-2}
\|D^{\gamma''}h\|_{L_2(B_{r_2})} \right),
$$ 
where $r < r_1 < r_2 < R$. We repeat this argument as
many times as we need. The lemma is proved.
\end{proof}

Denote by $h_{x}$ a generic derivative $h_{x^j}$, $j = 1,\dots,d$, and
$h_{x'}$ a generic derivative $h_{x^j}$, $j = 2,\dots,d$. 
Thus, for example, $h_{xx'}$
can be $h_{x^jx^k}$ where $j \in \{1,2,\dots,d\}$ and $k \in
\{2,\dots,d\}$.

\begin{lem} \label{sup esti} Let $h$ be a sufficiently smooth function
$h$ defined on $B_4$ such that $Lh = 0$ in $B_4$. Then we have 
$$
\sup_{B_1} |h_{xxx'}|^2
\le N \int_{B_3} |h|^2 \, dx,
$$ 
where  $N = N(d, \delta)$.
\end{lem}

\begin{proof} 
Imagine that we have
\begin{equation}                                   
                                                \label{maincase01}
\sup_{B_1}|h_{xx}| \le  N(d, \delta) 
\| h \|_{L_2(B_{5/2})}.
\end{equation} 
Then using   the fact that
$Lh_{x'} = 0$ we would obtain
$$
\sup_{B_1}|h_{x'xx}| \le N \| h_{x'} \|_{L_2(B_{5/2})} 
$$ 
and it would only remain to appeal to  Lemma \ref{conse esti}.

Therefore, it suffices to prove \eqref{maincase01}.
To do that, we first fix an integer $k$
such that $k - (d-1)/2 > 0$. Then due to the Sobolev embedding theorem,
we can find a constant $N$ such that, for each $-1 \le x^1 \le 1$,
$$
\sup_{|x'| \le 1} |h_{x'x^1x^1}(x^1, x')| \le N \|h_{x'x^1x^1}(x^1,
\cdot)\|_{W_2^k(B_1')}
$$ 
and
$$
\sup_{|x'| \le 1} |h_{x'x^1}(x^1, x')| \le N \|h_{x'x^1}(x^1,
\cdot)\|_{W_2^k(B_1')},
$$ 
where $B_1' = \{ x' \in \bR^{d-1} : |x'| \le 1 \}$. Set
$g $ to be either $h_{x'x^1x^1} $ or $h_{x'x^1} $.
Then 
$$
\begin{aligned}
\int_{-1}^{1} \sup_{|x'| \le 1} | g(x^1, x') |^2 \, dx^1
\phantom{i}&\le N \int_{-1}^{1}  \|g(x^1, \cdot)\|^2_{W_2^k(B_1')} \,
dx^1 \\
\phantom{i}&\le N \sum_{\substack{|\gamma| \le k+3 \\ 1 \le \gamma^1 \le
2}}
\| D^{\gamma} h \|^2_{L_2(B_2)}.
\end{aligned}
$$ 
 From this and Lemma \ref{conse esti} we have
\begin{equation}\label{obser01}
\int_{-1}^{1} \sup_{|x'| \le 1} | h_{x'x^1x^1} |^2 \, dx^1 +
\int_{-1}^{1} \sup_{|x'| \le 1} | h_{x'x^1} |^2 \, dx^1
\le N \| h \|^2_{L_2(B_{5/2})},
\end{equation} where $N$ depends only on  $d$ and $\delta$. Now we
notice that, for $x^1, y^1 \in [-1,1]$,
\begin{multline*}
 \sup_{|x'| \le 1}|h_{x'x^1}(x^1,x')| - \sup_{|x'| \le
1}|h_{x'x^1}(y^1,x')|   \\
\le \sup_{|x'| \le 1} \left| h_{x'x^1}(x^1,x') - h_{x'x^1}(y^1,x')
\right| 
\le \int_{x^1}^{y^1} \sup_{|x'| \le 1} |h_{x' x^1 x^1}(t, x')| \, dt \\
\le |x^1 - y^1|^{1/2} \left( \int_{-1}^{1} \sup_{|x'| \le 1} |h_{x' x^1
x^1}(t, x')|^2 \, dt \right)^{1/2}.
\end{multline*} This and \eqref{obser01} imply
$$
\sup_{|x'| \le 1}|h_{x'x^1}(x^1,x')|
\le N \| h \|_{L_2(B_{5/2})} \, |x^1 - y^1|^{1/2} + \sup_{|x'| \le
1}|h_{x'x^1}(y^1,x')|.
$$ 
Take integrals of both sides with respect to $y^1$, and
take a supremum over $x^1$. Then
$$
\begin{aligned}
\sup_{x \in B_1}|h_{x'x^1}(x)|
\phantom{i}& \le N \| h \|_{L_2(B_{5/2})} + \int_{-1}^{1} \sup_{|x'| \le
1}|h_{x'x^1}(y^1,x')| \, dy^1 \\
\phantom{i}& \le N \| h \|_{L_2(B_{5/2})},
\end{aligned}
$$ 
where the last inequality follows from \eqref{obser01},
and 
$N$ depends only on  $d$ and $\delta$. Similarly, we follows the
same steps as above with $h_{x'x'x^1}$ and $h_{x'x'}$ in place of $h_{x'
x^1 x^1} $ and $h_{x'x^1} $, respectively. Therefore, we have 
$$
\sup_{x \in B_1}|h_{x' x}(x)| \le  N(d,\delta) \| h \|_{L_2(B_{5/2})} .
$$ 
Finally, using the fact that $a^{11}h_{x^1x^1} = -\sum_{j
\ne 1 \, \text{or} \, k \ne 1} a^{jk}h_{x^jx^k}$, we finish the proof of
\eqref{maincase01}.
\end{proof}

Denote by $(u)_{B_r(x_0)}$ the average value of 
a function $u$ over $B_r(x_0)$,
that is,
$$
 (u)_{B_r(x_0)} = \eint_{B_r(x_0)} u(x) \,
dx = \frac{1}{|B_r|}\int_{B_r(x_0)} u(x) \, dx.
$$

Let $u \in C_0^{\infty}(\bR^d)$ and $f:= Lu$.
Assume that
$a^{jk}(x^1)$ are infinitely differentiable as functions of $x^1 \in \bR$.
Then we can find a sufficiently smooth function $h$ defined on $B_4$ such
that
$$
\left\{ \begin{aligned} Lh \phantom{i}&= 0 \quad \text{in} \quad B_4\\ h
\phantom{i}&= u \quad \text{on} \quad \partial B_4
\end{aligned} \right..
$$ 
For this solution $h$, we   establish the following inequality.

\begin{lem}                                   \label{sup esti_01}  
There exists a constant
 $N = N(d, \delta)$  such that 
$$
\sup_{B_1} |h_{xxx'}|^{2}
\le N   \int_{B_4} |f|^{2 } \, dx   + N
 \int_{B_4} |u_{xx}|^{2  } \, dx  .
$$
\end{lem}

\begin{proof} Define
$$
\begin{aligned}
\tilde{u} := u - u_{B_4} - (u_{x^i})_{B_4} x^i
\quad \text{in} \quad B_4, \\
\tilde{h} := h - u_{B_4} - (u_{x^i})_{B_4} x^i
\quad \text{in} \quad B_4.
\end{aligned}
$$ 
Then
$$
L \tilde{u} = f, \quad
L \tilde{h} = 0 \quad \text{in} \quad B_4
\quad \text{and} \quad
 \tilde{h} = \tilde{u}  \quad \text{on} \quad \partial B_4.
$$
 By Lemma \ref{sup esti} we see that
$$
\sup_{B_1}|h_{xxx'}|^2 = \sup_{B_1} |\tilde{h}_{xxx'}|^2
\le N \int_{B_3} |\tilde{h}|^2 \, dx.
$$

Let $\eta$ be a function in $C_0^{\infty}(\bR^d)$ such that
$\eta(x) = 0$ in $B_3$ and $\eta(x) = 1$ at $\partial B_4$. Then
$\tilde{h} - \eta \tilde{u} \in W_2^2(B_4)$ and $\tilde{h} - \eta
\tilde{u} = 0$ on $\partial B_4$. Therefore, by 
Lemma \ref{L_2 estimate}
$$
\int_{B_3} |\tilde{h}|^2 \, dx = \int_{B_3} |\tilde{h} - \eta
\tilde{u}|^2 \, dx
\le  N(d,\delta) 
 \int_{B_4} |L(\eta \tilde{u})|^{2 } \, dx  .
$$
 Note that 
$$
\begin{aligned} L(\eta \tilde{u})
\phantom{i}&= \eta L u + 2 a^{ij} \eta_{x^i} \tilde{u}_{x^j} + \tilde{u}
L \eta \\
\phantom{i}&= \eta f + 2 a^{ij} \eta_{x^i} (u_{x^j} - (u_{x^j})_{B_4}) +
(u - u_{B_4} - (u_{x^i})_{B_4} x^i) L \eta.
\end{aligned}
$$ 
Hence we have
$$
\int_{B_4} |L(\eta \tilde{u})|^{2 } \, dx 
 \le N \int_{B_4}( |f|^{2 } + |u_{x^j} -
(u_{x^j})_{B_4}|^{2  } )\,dx
$$
$$
+  N \int_{B_4}|u - u_{B_4} - (u_{x^i})_{B_4} x^i|^{2} \, dx  
 \le N \int_{B_4} |f|^{2  } \, dx +  N \int_{B_4}
|u_{xx}|^{2  } \, dx,
$$
 where the last inequality follows from Lemmas 3.1 and 3.2
in \cite{Krylov_2005}, and $N$ depends only on  $d$ and $\delta$.
The lemma is proved.
\end{proof}

\begin{lem}\label{sharp h} Let   $\kappa \ge 4$, and $r > 0$. 
Also let $a^{jk}(x^1)$ be infinitely differentiable. For a given $u \in
C_0^{\infty}(\bR^d)$, we find a smooth function $h$ defined
on $B_{\kappa r}$ such that $Lh = 0$ in $B_{\kappa r}$ and
$h = u$ on $\partial B_{\kappa r}$. Then there exists a constant 
 $N = N(d, \delta)$  such that
\begin{equation}\label{h inequal}
\eint_{B_r} |h_{xx'} - \left(h_{xx'}\right)_{B_r}|^2 \, dx 
\le N \kappa^{-2} \left[ \left(|Lu|^{2  }\right)_{B_{\kappa
r}}  + \left(|u_{xx}|^{2  }\right)_{B_{\kappa
r}}  \right].
\end{equation}
\end{lem}

\begin{proof} Using the dilation argument as in the proof of
Lemma
\ref{L_2 estimate}, we see that we need to prove only the case $r = 1$.
In that case we first observe that by using the same dilation
argument and Lemma \ref{sup esti_01}, we have
$$
\sup_{B_{\kappa/4}} | h_{xxx'} |^2 
\le N \kappa^{-2} \left[ \left( |Lu|^{2  }
\right)_{B_{\kappa}}  + 
\left( |u_{xx}|^{2  } \right)_{B_{\kappa}}  \right],
$$
 where $N$ depends only on  $d$ and $\delta$. 
Now we need only observe that $\kappa/4 \ge 1$,
$r=1$, and the left
hand side of the inequality \eqref{h inequal} is not greater than a
constant times $\sup_{B_1} | h_{xxx'} |^2$. The lemma is proved.
\end{proof}

Using the results obtained above, we will finally arrive at

\begin{lem}                                           \label{f_sh_esti} 
There exists a
constant  $N = N(d,\delta)$  such that, for any $\kappa \ge 4$, $r
> 0$, and $u \in C_0^{\infty}(\bR^d)$, we have
\begin{equation}         
                                                 \label{estimate_sh}
\eint_{B_r} |u_{xx'} - 
\left(u_{xx'}\right)_{B_r}|^2 \, dx 
\le  N \kappa^{d } \left( |Lu|^2 \right)_{B_{\kappa r}}   + N
\kappa^{-2}  
\left( |u_{xx}|^{2  } \right)_{B_{\kappa r}}  .
\end{equation}
\end{lem}

\begin{proof} We can assume that $a^{jk}(x^1)$ are infinitely
differentiable. In that case, we find a sufficiently smooth $h$ defined
on $B_{\kappa r}$ such that
$Lh = 0$ in $B_{\kappa r}$ and $h = u$ on $\partial B_{\kappa r}$. Note
that 
$L(u-h) = Lu$ in $B_{\kappa r}$ and $u-h = 0$ on $\partial B_{\kappa r}$.
From  Lemma \ref{sharp h} we have
\begin{equation}\label{cal 001}
\eint_{B_r} |h_{xx'} - \left(h_{xx'}\right)_{B_r}|^2 \, dx 
\le N \kappa^{-2} \left[ \left( |Lu|^{2  } \right)_{B_{\kappa
r}}  + 
\left( |u_{xx}|^{2  } \right)_{B_{\kappa r}}  \right].
\end{equation} On the other hand, from   estimate \eqref{case lambda0}
we have
$$
\int_{B_r} |u_{xx'} - h_{xx'}|^2 \, dx
\le N \left( \int_{B_{\kappa r}} | Lu |^2 \, dx + r^{-4}(\kappa-1)^{-4}
\int_{B_{\kappa r}} |u-h|^2 \, dx \right).
$$ 
Moreover, by Lemma \ref{L_2 estimate} 
$$
\int_{B_{\kappa r}} |u - h|^2 \, dx 
\le N \, (\kappa r)^{4 }
 \int_{B_{\kappa r}} |Lu|^{2 } \, dx   .
$$
Hence
$$
\eint_{B_r} |u_{xx'} - h_{xx'}|^2 \, dx
\le N \kappa^{d} \left( |Lu|^2 \right)_{B_{\kappa r}} .
$$
 This and \eqref{cal 001} prove the inequality
\eqref{estimate_sh} with $\left(h_{xx'}\right)_{B_r}$ in place of
$\left(u_{xx'}\right)_{B_r}$. Now we need only notice that
$$
\eint_{B_r} |u_{xx'} - \left(u_{xx'}\right)_{B_r} |^2 \, dx \le  
\eint_{B_r} |u_{xx'} - \left(h_{xx'}\right)_{B_r} |^2 \, dx.
$$ 
The lemma is proved.
\end{proof}

\mysection{Proof of Theorem \ref{main_theorem}}
                                                       \label{W_p^2_case}

In this section we suppose that all
assumption of Theorem \ref{main_theorem}  are satisfied. Recall that
$$ 
Lu(x) = a^{jk}(x) u_{x^j x^k}(x) + b^{j}(x) u_{x^j}(x)
+ c(x) u(x),
$$
$$ 
L_0u(x) = a^{jk}(x) u_{x^j x^k}(x).
$$ 
We use the maximal and sharp functions 
 given by
$$ 
M g (x) = \sup_{r>0} \eint_{B_r(x)} |g(y)| \, dy,
$$
$$
g^{\#}(x) = \sup_{r>0} \eint_{B_r(x)} |g(y) -
(g)_{B_r(x)}| \, dy.
$$

\begin{thm} Let  $\mu$, $\nu \in (1,\infty)$, $1/\mu + 1/\nu =
1$, and $R \in (0, \infty)$. Then there exists a constant 
 $N = N(d,\delta,\mu)$  such that, for any $u \in C_0^{\infty}(B_R)$,
we have
 
\begin{multline}
                                                \label{sharp_esti}
\left(u_{xx'}\right)^{\#}
\le  N \left(a^{\#(x')}_R
\right)^{\alpha} \left[
M(|u_{xx}|^{2  \mu}) \right]^{\beta} \\ + N \left[
M(|L_0u|^{2  }) \right]^{1/(d+2) } 
\left[ M(|u_{xx}|^{2  }) \right]^{d/ (2d+4)},
\end{multline}
where $\alpha=\nu^{-1}(d+2)^{-1}$, $\beta=2^{-1}\mu^{-1}$.
 
\end{thm}

\begin{proof} Fix $\kappa \ge 4$, $r \in (0, \infty)$, and $x_0 = (x_0^1,
x_0') \in \bR^d$. Introduce 
$$
\bar{a}^{jk}(x^1) = \frac{1}{|B'_{\kappa r}|} \int_{B'_{\kappa r}(x'_0)}
a^{jk}(x^1,y') \, dy' \quad \text{if} \quad \kappa r < R,
$$
$$
\bar{a}^{jk}(x^1) = \frac{1}{|B'_R|} \int_{B'_R} a^{jk}(x^1,y') \, dy'
\quad \text{if} \quad \kappa r \ge R,
$$
$$
\mathcal{A} =   M(|L_0 u|^{2 })(x_0)  ,
\qquad
\mathcal{B} =  M(|u_{xx}|^{2  })(x_0)  ,
$$
$$
\mathcal{C} = \left( M(|u_{xx}|^{2   \mu})(x_0) \right)^{1/ 
\mu}.
$$

Set $\bar{L}_0 u = \bar{a}^{jk}(x^1) u_{x^j x^k}$.  Then  Lemma
\ref{f_sh_esti} along with  the
fact that $\kappa \ge 4$ allows us to obtain
\begin{equation}                                 \label{sim_ineq_01}
\begin{aligned}
\phantom{i}&\eint_{B_r(x_0)} |u_{xx'} -
\left(u_{xx'}\right)_{B_r(x_0)}|^2 \, dx \\
\phantom{i}&\le  N \kappa^{d } \left( |\bar{L}_0 u|^{2  }
\right)_{B_{\kappa r}(x_0)}   + N \kappa^{-2} \left(|u_{xx}|^{2
\ } \right)_{B_{\kappa r}(x_0)} 
\end{aligned}
\end{equation} 
for $\kappa \ge 4$, where $N$ depends only on 
 $d$ and $\delta$. Note that
\begin{equation}                                       \label{diff01}
\int_{B_{\kappa r}(x_0)}
 | \bar{L}_0 u |^{2  } \, dx
\le 2\int_{B_{\kappa r}(x_0)} | \bar{L}_0 u - L_0 u
|^{2 } \, dx + 2  \int_{B_{\kappa r}(x_0)} |L_0
u|^{2} \, dx
\end{equation} 
and
$$
\int_{B_{\kappa r}(x_0)} | \bar{L}_0 u - L_0 u |^{2 } \, dx =
\int_{B_{\kappa r}(x_0) \cap B_R} | \bar{L}_0 u - L_0 u |^{2 } \, dx
$$
$$
\le \left(\int_{B_{\kappa r}(x_0) \cap B_R} |\bar{a}   -
a |^{2  \nu} \, dx \right)^{1/\nu}
\left(\int_{B_{\kappa r}(x_0)} |u_{x x } |^{2  \mu} \, dx
\right)^{1/\mu}  := I^{1/\nu} J^{1/\mu}.
$$
  If $\kappa r < R$, we have
\begin{multline*} 
I \le N \int_{x_0^1 - \kappa r}^{x_0^1 + 
\kappa r} \int_{B'_{\kappa
r}(x_0')} |\bar{a} (x^1) - a (x^1,x')| \, dx' \, dx^1 \\
\le N (\kappa r)^d a^{\#(x')}_{\kappa r} \le N (\kappa r)^d
a^{\#(x')}_R.
\end{multline*} 
In case $\kappa r \ge R$ 
\begin{multline*} I 
\le N \int_{-R}^{R} \int_{B'_R}  |\bar{a}(x^1) - a(x^1,x')| 
\, dx' \, dx^1 \\
\le N R^d a^{\#(x')}_R \le N (\kappa r)^d a^{\#(x')}_R.
\end{multline*} Hence
$$
\int_{B_{\kappa r}(x_0)} | \bar{L}_0 u - L_0 u |^{2 } \, dx
\le N (\kappa r)^{d/\nu} \left( a^{\#(x')}_R \right)^{1/\nu}
\left( \int_{B_{\kappa r}(x_0)} |u_{xx}|^{2  \mu} \, dx
\right)^{1/\mu}.
$$ 
From  this and \eqref{diff01} it follows that
$$
\left( |\bar{L}_0 u|^{2  } \right)_{B_{\kappa r}(x_0)} 
\le N \left[ \left(a^{\#(x')}_R\right)^{1/  \nu} 
\left( |u_{xx}|^{2   \mu} \right)_{B_{\kappa r}(x_0)}^{1/ \mu}
+ \left( |L_0 u|^{2  } \right)_{B_{\kappa r}(x_0)} 
\right].
$$
This and \eqref{sim_ineq_01} allow us to have
\begin{multline*}
\eint_{B_r(x_0)} |u_{xx'} - \left(u_{xx'}\right)_{B_r(x_0)}|^2 \, dx 
\le  N \kappa^{d } \left(a^{\#(x')}_R\right)^{1/ \nu} 
\left( |u_{xx}|^{2  \mu} \right)_{B_{\kappa r}(x_0)}^{1/  \mu}
\\ + N \kappa^{d } \left( |L_0 u|^{2  } \right)_{B_{\kappa
r}(x_0)}   + N \kappa^{-2} \left(|u_{xx}|^{2  }
\right)_{B_{\kappa r}(x_0)}  \\
\le N \kappa^{d } \left(a^{\#(x')}_R\right)^{1/  \nu} \mathcal{C} +
N \kappa^{d } \mathcal{A} + N \kappa^{-2} \mathcal{B},
\end{multline*}
for all $r>0$ and $\kappa \ge 4$. In addition, the above
inequality is also true for $0 < \kappa < 4$ since then
$$
\eint_{B_r(x_0)} |u_{xx'} - \left(u_{xx'}\right)_{B_r(x_0)}|^2 \, dx
\leq\eint_{B_r(x_0)} |u_{xx'}|^{2}\,dx\leq
\mathcal{B}\leq 16\kappa^{-2}\mathcal{B}.
$$  
  By taking the supremum with respect
to $r > 0$, and then minimizing with respect to $\kappa > 0$, we have
\begin{multline*}
\left[u^{\#}_{xx'}(x_0)\right]^2
\le N \left(\left(a^{\#(x')}_R\right)^{1/\nu} \mathcal{C} +
\mathcal{A} \right)^{\frac{2}{d+2}}
\mathcal{B}^{\frac{d}{d+2}} \\
\le N \left(a^{\#(x')}_R\right)^{\frac{2}{\nu (d+2)}} \,
\mathcal{C}^{\frac{2}{d+2}} \mathcal{B}^{\frac{d}{d+2}}  + N \,
\mathcal{A}^{\frac{2}{d+2}} \, \mathcal{B}^{\frac{d}{d+2}},
\end{multline*} 
where $N = N(d,\delta,\mu)$.  Notice that
$\mathcal{B} \le \mathcal{C}$. Thus by replacing $\mathcal{B}$ with
$\mathcal{C}$ in the first term on the right we finish the proof.
\end{proof}

\begin{cor} For $p > 2$, there exist constants 
 $R = R(d, \delta, p,\omega)$ and $N = N(d, \delta, p)$  such that, for any $u \in C_0^{\infty}(B_R)$, 
we have
\begin{equation}\label{local_esti}
\| u_{xx} \|_{L_p} \le N \| L_0 u \|_{L_p}.
\end{equation}
\end{cor}

\begin{proof}  Choose real numbers $\mu > 1$ such that $p
> 2 \mu$.  Then we use the inequality \eqref{sharp_esti} along with the
Fefferman-Stein theorem on sharp functions and the Hardy-Littlewood
maximal function theorem. We also use H\"{o}lder's inequality to have
(note that  $p/2 \mu > 1$ and $p/2 > 1$ )
\begin{equation}\label{esti_xx'}
\| u_{xx'} \|_{L_p}
\le N \left(a^{\#(x')}_R\right)^{\frac{1}{\nu (d+2)}}
\|u_{xx}\|_{L_p} + N \| L_0 u \|^{\frac{2}{d+2}}_{L_p}
\|u_{xx}\|^{\frac{d}{d+2}}_{L_p},
\end{equation}
 where $1/\mu + 1/\nu = 1$, and $N$ depends only on 
 $d$, $\delta$, and $p$. Since  
$$ 
u_{x^1 x^1} = L_0 u - \sum_{j \ne 1 \, \text{or} \, k \ne 1}
\frac{a^{jk}}{a^{11}} u_{x^j x^k},
$$ 
by using \eqref{esti_xx'} we arrive at
\begin{equation}
                                                       \label{12.6.1}
\| u_{xx} \|_{L_p}
\le N_{1} \left(a^{\#(x')}_R\right)^{\frac{1}{\nu (d+2)}}
\|u_{xx}\|_{L_p}  + N \| L_0 u \|_{L_p} + N \| L_0 u
\|^{\frac{2}{d+2}}_{L_p} \|u_{xx}\|^{\frac{d}{d+2}}_{L_p}.
\end{equation}
 We now invoke Assumption \ref{assum_02} by which we can
choose a sufficiently small $R$ such that
$$ 
N_{1} \left(a^{\#(x')}_R\right)^{\frac{1}{\nu
(d+2)}}  \le 1/2.
$$
 Then we have
$$
\frac{1}{2}\| u_{xx} \|_{L_p}
\le N \| L_0 u \|_{L_p} + N \| L_0 u \|^{\frac{2}{d+2}}_{L_p}
\|u_{xx}\|^{\frac{d}{d+2}}_{L_p},
$$
 which implies  \eqref{local_esti}.
\end{proof}

\begin{proof}[\textit{\textbf{Proof of Theorem ~\ref{main_theorem}}}]
Since we have an $L_p$-estimate for functions with small compact support,
we can just follow the standard argument, which can be found in
\cite{Krylov_2005}.
\end{proof}

\bibliographystyle{plain}

\def\cprime{$'$}\def\cprime{$'$} \def\cprime{$'$} \def\cprime{$'$}

\end{document}